\documentclass[a4paper,11pt]{amsart}
\usepackage{amsmath,amssymb,amsfonts,latexsym}

\newtheorem{theorem}{Theorem}[section]

\input{tcilatex}

\begin{document}
\title[\textbf{Extended} $q$\textbf{-Dedekind-type DC sums}]{\textbf{Extended%
} $q$\textbf{-Dedekind-type Daehee-Changhee sums\ associated with Extended }$%
q$\textbf{-Euler polynomials}}
\author[\textbf{S. Araci}]{\textbf{Serkan Araci}}
\address{\textbf{University of Gaziantep, Faculty of Science and Arts,
Department of Mathematics, 27310 Gaziantep, TURKEY}}
\email{\textbf{mtsrkn@hotmail.com; saraci88@yahoo.com.tr; mtsrkn@gmail.com}}
\author[\textbf{M. Acikgoz}]{\textbf{Mehmet Acikgoz}}
\address{\textbf{University of Gaziantep, Faculty of Science and Arts,
Department of Mathematics, 27310 Gaziantep, TURKEY}}
\email{\textbf{acikgoz@gantep.edu.tr}}

\begin{abstract}
In the present paper, our goal is to introduce a $p$-adic continuous
function for an odd prime to inside a $p$-adic $q$-analogue of the higher
order Dedekind-type sums with weight $\alpha $ related to Extended $q$-Euler
polynomials by using $p$-adic $q$-integral in the $p$-adic integer ring.

\vspace{2mm}\noindent \textsc{2010 Mathematics Subject Classification.}
11S80, 11B68.

\vspace{2mm}

\end{abstract}

\thanks{}
\maketitle




\section{\textbf{Introduction}}


Assume that $p$ be a fixed odd prime number. Throughout this paper $%
\mathbb{Z}
_{p}$, $%
\mathbb{Q}
_{p}$, $%
\mathbb{C}
$ and $%
\mathbb{C}
_{p}$ will, respectively, denote the ring of $p$-adic rational integers, the
field of $p$-adic rational numbers, the complex numbers and the completion
of algebraic closure of $%
\mathbb{Q}
_{p}$.

Let $v_{p}$ be normalized exponential valuation of $%
\mathbb{C}
_{p}$ by%
\begin{equation*}
\left\vert p\right\vert _{p}=p^{-v_{p}\left( p\right) }=\frac{1}{p}\text{.}
\end{equation*}

When one talks of $q$-extension, $q$ is variously considered as an
indeterminate, a complex number $q\in 
\mathbb{C}
$ or $p$-adic number $q\in 
\mathbb{C}
_{p}$. If $q\in 
\mathbb{C}
$, we assume that $\left\vert q\right\vert <1$. If $q\in 
\mathbb{C}
_{p}$, we assume that $\left\vert 1-q\right\vert _{p}<1$ (see, [1-15]).

A $q$-analogue of $p$-adic Haar distribution is defined by Kim as follows:
for any postive integer $n$,%
\begin{equation*}
\mu _{q}\left( a+p^{n}%
\mathbb{Z}
_{p}\right) =\left( -q\right) ^{a}\frac{\left( 1+q\right) }{1+q^{p^{n}}}
\end{equation*}%
for $0\leq a<p^{n}$ and this can be extended to a measure on $%
\mathbb{Z}
_{p}$ (for details, see [1-7]).

Extended $q$-Euler polynomials are defined as%
\begin{equation}
\widetilde{E}_{n,q}^{\left( \alpha \right) }\left( x\right) =\int_{%
\mathbb{Z}
_{p}}\left( \frac{1-q^{\alpha \left( x+\xi \right) }}{1-q^{\alpha }}\right)
^{n}d\mu _{q}\left( \xi \right)  \label{equation 1}
\end{equation}%
for $n\in 
\mathbb{Z}
_{+}:=\left\{ 0,1,2,3,\cdots \right\} $. We note that 
\begin{equation*}
\lim_{q\rightarrow 1}\widetilde{E}_{n,q}^{\left( \alpha \right) }\left(
x\right) =E_{n}\left( x\right)
\end{equation*}%
where $E_{n}\left( x\right) $ are $n$-th Euler polynomials, which are
defined by the rule:%
\begin{equation*}
\sum_{n=0}^{\infty }E_{n}\left( x\right) \frac{t^{n}}{n!}=e^{tx}\frac{2}{%
e^{t}+1},\text{ }\left\vert t\right\vert <\pi
\end{equation*}%
(for details, see \cite{Ryoo}). Taking $x=0$ into (\ref{equation 1}), then
we have $\widetilde{E}_{n,q}^{\left( \alpha \right) }\left( 0\right) :=%
\widetilde{E}_{n,q}^{\left( \alpha \right) }$ are called extended $q$-Euler
numbers.

Extended $q$-Euler numbers and polynomials have the following equalities%
\begin{eqnarray}
\widetilde{E}_{n,q}^{\left( \alpha \right) } &=&\frac{1+q}{\left(
1-q^{\alpha }\right) ^{n}}\sum_{l=0}^{n}\binom{n}{l}\left( -1\right) ^{l}%
\frac{1}{1+q^{\alpha l+1}}\text{,}  \label{equation 2} \\
\widetilde{E}_{n,q}^{\left( \alpha \right) }\left( x\right) &=&\frac{1+q}{%
\left( 1-q^{\alpha }\right) ^{n}}\sum_{l=0}^{n}\binom{n}{l}\left( -1\right)
^{l}\frac{q^{\alpha lx}}{1+q^{\alpha l+1}}\text{,}  \label{equation 3} \\
\widetilde{E}_{n,q}^{\left( \alpha \right) }\left( x\right) &=&\sum_{l=0}^{n}%
\binom{n}{l}q^{\alpha lx}\widetilde{E}_{l,q}^{\left( \alpha \right) }\left( 
\frac{1-q^{\alpha x}}{1-q^{\alpha }}\right) ^{n-l}\text{.}
\label{equation 4}
\end{eqnarray}

Also, for $d\in 
\mathbb{N}
$ with $d\equiv 1\left( \func{mod}2\right) $ 
\begin{equation}
\widetilde{E}_{n,q}^{\left( \alpha \right) }\left( x\right) =\left( \frac{1+q%
}{1+q^{d}}\right) \left( \frac{1-q^{\alpha d}}{1-q^{\alpha }}\right)
^{n}\sum_{a=0}^{d-1}\left( -1\right) ^{a}\widetilde{E}_{n,q}^{\left( \alpha
\right) }\left( \frac{x+a}{d}\right) \text{,}  \label{equation 5}
\end{equation}

(for more information, see \cite{Ryoo}).

For any positive integer $h,k$ and $m$, Dedekind-type DC sums are given by
Kim in \cite{Kim 1}, \cite{Kim 2} and \cite{Kim 3} as follows:%
\begin{equation*}
S_{m}\left( h,k\right) =\sum_{M=1}^{k-1}\left( -1\right) ^{M-1}\frac{M}{k}%
\overline{E}_{m}\left( \frac{hM}{k}\right)
\end{equation*}%
where $\overline{E}_{m}\left( x\right) $ are the $m$-th periodic Euler
function.

Recently, weighted $q$-Bernoulli numbers and polynomials was firstly defined
by Taekyun Kim in \cite{Kim 8}. After, many mathematicians, by utilizing
from Kim's paper \cite{Kim 8}, have introduced a new concept in Analytic
numbers theory as weighted $q$-Bernoulli, weighted $q$-Euler, weighted $q$%
-Genocchi polynomials in \cite{Rim}, \cite{Ryoo}, \cite{Araci 1} and \cite%
{Araci 3}. Also, Kim derived some interesting properties for Dedekind-type
DC sums. He considered a $p$-adic continuous function for an odd prime
number to contain a $p$-adic $q$-analogue of the higher order Dedekind-type
DC sums $k^{m}S_{m+1}\left( h,k\right) $ in \cite{Kim 2}. In \cite{Simsek},
Simsek also studied to $q$-analogue of Dedekind-type sums. He also derived
their interesting properties.

By the same motivation, we, by using $p$-adic $q$-integral on $%
\mathbb{Z}
_{p}$, will construct weighted $p$-adic $q$-analogue of the higher order
Dedekind-type DC sums $k^{m}S_{m+1}\left( h,k\right) $.

\section{\textbf{Extended }$q$\textbf{-Dedekind-type Sums associated with
Extended }$q$\textbf{-Euler polynomials}}

Let $w$ be the $Teichm\ddot{u}ller$ character ($\func{mod}p$). For $x\in 
\mathbb{Z}
_{p}^{\ast }$ $:=%
\mathbb{Z}
_{p}/p%
\mathbb{Z}
_{p}$, set%
\begin{equation*}
\left\langle x:q\right\rangle =w^{-1}\left( x\right) \left( \frac{1-q^{x}}{%
1-q}\right) \text{.}
\end{equation*}

Let $a$ and $N$ be positive integers with $\left( p,a\right) =1$ and $p\mid
N $. We now consider the following%
\begin{equation*}
\widetilde{C}_{q}^{\left( \alpha \right) }\left( s,a,N:q^{N}\right)
=w^{-1}\left( a\right) \left\langle x:q^{\alpha }\right\rangle
^{s}\sum_{j=0}^{\infty }\binom{s}{j}q^{\alpha aj}\left( \frac{1-q^{\alpha N}%
}{1-q^{\alpha a}}\right) ^{j}\widetilde{E}_{j,q^{N}}^{\left( \alpha \right) }%
\text{.}
\end{equation*}

In particular, if $m+1\equiv 0(\func{mod}p-1)$, then%
\begin{eqnarray*}
\widetilde{C}_{q}^{\left( \alpha \right) }\left( m,a,N:q^{N}\right)
&=&\left( \frac{1-q^{\alpha a}}{1-q^{\alpha }}\right) ^{m}\sum_{j=0}^{m}%
\binom{m}{j}q^{\alpha aj}\widetilde{E}_{j,q^{N}}^{\left( \alpha \right)
}\left( \frac{1-q^{\alpha N}}{1-q^{\alpha a}}\right) ^{j} \\
&=&\left( \frac{1-q^{\alpha N}}{1-q^{\alpha }}\right) ^{m}\int_{%
\mathbb{Z}
_{p}}\left( \frac{1-q^{\alpha N\left( \xi +\frac{a}{N}\right) }}{1-q^{\alpha
N}}\right) ^{m}d\mu _{q^{N}}\left( \xi \right) \text{.}
\end{eqnarray*}

Thus, $\widetilde{C}_{q}^{\left( \alpha \right) }\left( m,a,N:q^{N}\right) $
is a continuous $p$-adic extension of 
\begin{equation*}
\left( \frac{1-q^{\alpha N}}{1-q^{\alpha }}\right) ^{m}\widetilde{E}%
_{m,q^{N}}^{\left( \alpha \right) }\left( \frac{a}{N}\right) \text{.}
\end{equation*}

Let $\left[ .\right] $ be the Gauss' symbol and let $\left\{ x\right\} =x-%
\left[ x\right] $. Thus, we are now ready to introduce $q$-analogue of the
higher order Dedekind-type DC sums $\widetilde{J}_{m,q}^{\left( \alpha
\right) }\left( h,k:q^{l}\right) $ by the rule: 
\begin{equation*}
\widetilde{J}_{m,q}^{\left( \alpha \right) }\left( h,k:q^{l}\right)
=\sum_{M=1}^{k-1}\left( -1\right) ^{M-1}\left( \frac{1-q^{\alpha M}}{%
1-q^{\alpha k}}\right) \int_{%
\mathbb{Z}
_{p}}\left( \frac{1-q^{\alpha \left( l\xi +l\left\{ \frac{hM}{k}\right\}
\right) }}{1-q^{\alpha l}}\right) ^{m}d\mu _{q^{l}}\left( \xi \right) \text{.%
}
\end{equation*}

If $m+1\equiv 0\left( \func{mod}p-1\right) $%
\begin{eqnarray*}
&&\left( \frac{1-q^{\alpha k}}{1-q^{\alpha }}\right)
^{m+1}\sum_{M=1}^{k-1}\left( -1\right) ^{M-1}\left( \frac{1-q^{\alpha M}}{%
1-q^{\alpha k}}\right) \int_{%
\mathbb{Z}
_{p}}\left( \frac{1-q^{\alpha k\left( \xi +\frac{hM}{k}\right) }}{%
1-q^{\alpha k}}\right) ^{m}d\mu _{q^{k}}\left( \xi \right) \\
&=&\sum_{M=1}^{k-1}\left( -1\right) ^{M-1}\left( \frac{1-q^{\alpha M}}{%
1-q^{\alpha }}\right) \left( \frac{1-q^{\alpha k}}{1-q^{\alpha }}\right)
^{m}\int_{%
\mathbb{Z}
_{p}}\left( \frac{1-q^{\alpha k\left( \xi +\frac{hM}{k}\right) }}{%
1-q^{\alpha k}}\right) ^{m}d\mu _{q^{k}}\left( \xi \right)
\end{eqnarray*}%
where $p\mid k$, $\left( hM,p\right) =1$ for each $M$. Thanks to (\ref%
{equation 1}), we easily state the following 
\begin{gather}
\left( \frac{1-q^{\alpha k}}{1-q^{\alpha }}\right) ^{m+1}\widetilde{J}%
_{m,q}^{\left( \alpha \right) }\left( h,k:q^{k}\right)  \label{equation 6} \\
=\sum_{M=1}^{k-1}\left( \frac{1-q^{\alpha M}}{1-q^{\alpha }}\right) \left( 
\frac{1-q^{\alpha k}}{1-q^{\alpha }}\right) ^{m}\left( -1\right) ^{M-1}\int_{%
\mathbb{Z}
_{p}}\left( \frac{1-q^{\alpha k\left( \xi +\frac{hM}{k}\right) }}{%
1-q^{\alpha k}}\right) ^{m}d\mu _{q^{k}}\left( \xi \right)  \notag \\
=\sum_{M=1}^{k-1}\left( -1\right) ^{M-1}\left( \frac{1-q^{\alpha M}}{%
1-q^{\alpha }}\right) \widetilde{C}_{q}^{\left( \alpha \right) }\left(
m,\left( hM\right) _{k}:q^{k}\right)  \notag
\end{gather}%
where $(hM)_{k}$ denotes the integer $x$ such that $0\leq x<n$ and $x\equiv
\alpha \left( \func{mod}k\right) $.

It is not difficult to indicate the following 
\begin{gather}
\int_{%
\mathbb{Z}
_{p}}\left( \frac{1-q^{\alpha \left( x+\xi \right) }}{1-q^{\alpha }}\right)
^{k}d\mu _{q}\left( \xi \right)  \label{equation 7} \\
=\left( \frac{1-q^{\alpha m}}{1-q^{\alpha }}\right) ^{k}\frac{1+q}{1+q^{m}}%
\sum_{i=0}^{m-1}\left( -1\right) ^{i}\int_{%
\mathbb{Z}
_{p}}\left( \frac{1-q^{\alpha m\left( \xi +\frac{x+i}{m}\right) }}{%
1-q^{\alpha m}}\right) ^{k}d\mu _{q^{m}}\left( \xi \right) \text{.}  \notag
\end{gather}

On account of (\ref{equation 6}) and (\ref{equation 7}), we easily see that%
\begin{gather}
\left( \frac{1-q^{\alpha N}}{1-q^{\alpha }}\right) ^{m}\int_{%
\mathbb{Z}
_{p}}\left( \frac{1-q^{\alpha N\left( \xi +\frac{a}{N}\right) }}{1-q^{\alpha
N}}\right) ^{m}d\mu _{q^{N}}\left( \xi \right)   \label{equation 8} \\
=\frac{1+q^{N}}{1+q^{Np}}\sum_{i=0}^{p-1}\left( -1\right) ^{i}\left( \frac{%
1-q^{\alpha Np}}{1-q^{\alpha }}\right) ^{m}\int_{%
\mathbb{Z}
_{p}}\left( \frac{1-q^{\alpha pN\left( \xi +\frac{a+iN}{pN}\right) }}{%
1-q^{\alpha pN}}\right) ^{m}d\mu _{q^{pN}}\left( \xi \right) \text{.}  \notag
\end{gather}

Because of (\ref{equation 6}), (\ref{equation 7}) and (\ref{equation 8}), we
develop the $p$-adic integration as follows: 
\begin{equation*}
\widetilde{C}_{q}^{\left( \alpha \right) }\left( s,a,N:q^{N}\right) =\frac{%
1+q^{N}}{1+q^{Np}}\sum_{\underset{a+iN\neq 0(\func{mod}p)}{0\leq i\leq p-1}%
}\left( -1\right) ^{i}\widetilde{C}_{q}^{\left( \alpha \right) }\left(
s,\left( a+iN\right) _{pN},p^{N}:q^{pN}\right) \text{.}
\end{equation*}

So,%
\begin{gather*}
\widetilde{C}_{q}^{\left( \alpha \right) }\left( m,a,N:q^{N}\right) =\left( 
\frac{1-q^{\alpha N}}{1-q^{\alpha }}\right) ^{m}\int_{%
\mathbb{Z}
_{p}}\left( \frac{1-q^{\alpha N\left( \xi +\frac{a}{N}\right) }}{1-q^{\alpha
N}}\right) ^{m}d\mu _{q^{N}}\left( \xi \right) \\
-\left( \frac{1-q^{\alpha Np}}{1-q^{\alpha }}\right) ^{m}\int_{%
\mathbb{Z}
_{p}}\left( \frac{1-q^{\alpha pN\left( \xi +\frac{a+iN}{pN}\right) }}{%
1-q^{\alpha pN}}\right) ^{m}d\mu _{q^{pN}}\left( \xi \right)
\end{gather*}%
where $\left( p^{-1}a\right) _{N}$ denotes the integer $x$ with $0\leq x<N$, 
$px\equiv a\left( \func{mod}N\right) $ and $m$ is integer with $m+1\equiv 0(%
\func{mod}p-1)$. Therefore, we procure the following%
\begin{gather*}
\sum_{M=1}^{k-1}\left( -1\right) ^{M-1}\left( \frac{1-q^{\alpha M}}{%
1-q^{\alpha }}\right) \widetilde{C}_{q}^{\left( \alpha \right) }\left(
m,hM,k:q^{k}\right) \\
=\left( \frac{1-q^{\alpha k}}{1-q^{\alpha }}\right) ^{m+1}\widetilde{J}%
_{m,q}^{\left( \alpha \right) }\left( h,k:q^{k}\right) -\left( \frac{%
1-q^{\alpha k}}{1-q^{\alpha }}\right) ^{m+1}\left( \frac{1-q^{\alpha kp}}{%
1-q^{\alpha k}}\right) \widetilde{J}_{m,q}^{\left( \alpha \right) }\left(
\left( p^{-1}h\right) ,k:q^{pk}\right)
\end{gather*}%
where $p\nmid k$ and $p\nmid hm$ for each $M$. Thus, we give the following
definition, which seems interesting for further studying in theory of
Dedekind sums.

\begin{definition}
Let $h,k$ be positive integer with $\left( h,k\right) =1$, $p\nmid k$. For $%
s\in 
\mathbb{Z}
_{p},$ we define $p$-adic Dedekind-type DC sums as follows:  
\begin{equation*}
\widetilde{J}_{p,q}^{\left( \alpha \right) }\left( s:h,k:q^{k}\right)
=\sum_{M=1}^{k-1}\left( -1\right) ^{M-1}\left( \frac{1-q^{\alpha M}}{%
1-q^{\alpha }}\right) \widetilde{C}_{q}^{\left( \alpha \right) }\left(
m,hM,k:q^{k}\right) \text{.}
\end{equation*}
\end{definition}

As a result of the above definition, we state the following theorem.

\begin{theorem}
For $m+1\equiv 0(\func{mod}p-1)$ and $\left( p^{-1}a\right) _{N}$ denotes
the integer $x$ with $0\leq x<N$, $px\equiv a\left( \func{mod}N\right) $,
then, we have 
\begin{gather*}
\widetilde{J}_{p,q}^{\left( \alpha \right) }\left( s:h,k:q^{k}\right)
=\left( \frac{1-q^{\alpha k}}{1-q^{\alpha }}\right) ^{m+1}\widetilde{J}%
_{m,q}^{\left( \alpha \right) }\left( h,k:q^{k}\right) \\
-\left( \frac{1-q^{\alpha k}}{1-q^{\alpha }}\right) ^{m+1}\left( \frac{%
1-q^{\alpha kp}}{1-q^{\alpha k}}\right) \widetilde{J}_{m,q}^{\left( \alpha
\right) }\left( \left( p^{-1}h\right) ,k:q^{pk}\right) \text{.}
\end{gather*}
\end{theorem}

In the special case $\alpha =1$, our applications in theory of Dedekind sums
resemble Kim's results in \cite{Kim 2}. These results seem to be interesting
for further studies in \cite{Kim 1}, \cite{Kim 3} and \cite{Simsek}.



\end{document}